\documentclass[11pt,reqno]{amsart}

\usepackage[T1]{fontenc}
\usepackage{lmodern}
\usepackage{microtype}
\usepackage{amsmath,amssymb,amsthm,mathtools}
\usepackage{booktabs}
\usepackage{hyperref}

\hypersetup{
  pdftitle={Profinite and p-Adic Analytic Bruck Loops as Symmetric Sections of Congruence Subgroups},
  pdfauthor={Raffaello Caserta, PhD},
  pdfsubject={p-adic Lie groups, pro-p groups, Bruck loops, A\_l-loops, symmetric spaces},
  pdfkeywords={Bruck loop, A\_l-loop, Kikkawa loop, pro-p group, p-adic Lie group, Cayley transform, Catalan numbers, congruence subgroup, nilpotency class},
  colorlinks=true,
  linkcolor=blue,
  citecolor=blue,
  urlcolor=blue
}

\newtheorem{theorem}{Theorem}[section]
\newtheorem{proposition}[theorem]{Proposition}
\newtheorem{lemma}[theorem]{Lemma}

\theoremstyle{definition}

\theoremstyle{remark}
\newtheorem{remark}[theorem]{Remark}

\DeclareMathOperator{\Aut}{Aut}
\DeclareMathOperator{\Ad}{Ad}
\DeclareMathOperator{\GL}{GL}
\DeclareMathOperator{\SL}{SL}
\DeclareMathOperator{\PGL}{PGL}
\DeclareMathOperator{\PSL}{PSL}
\DeclareMathOperator{\SO}{SO}
\DeclareMathOperator{\SU}{SU}
\DeclareMathOperator{\OO}{O}
\DeclareMathOperator{\Sym}{Sym}
\DeclareMathOperator{\Stab}{Stab}
\DeclareMathOperator{\cl}{cl}
\DeclareMathOperator{\rk}{rk}
\DeclareMathOperator{\id}{id}
\DeclareMathOperator{\tr}{tr}

\newcommand{\Z}{\mathbb{Z}}
\newcommand{\Q}{\mathbb{Q}}
\newcommand{\F}{\mathbb{F}}
\newcommand{\R}{\mathbb{R}}
\newcommand{\Zp}{\mathbb{Z}_p}
\newcommand{\Qp}{\mathbb{Q}_p}
\newcommand{\Fp}{\mathbb{F}_p}
\newcommand{\Oint}{\mathcal{O}}
\newcommand{\h}{\mathfrak{h}}
\newcommand{\m}{\mathfrak{m}}
\newcommand{\so}{\mathfrak{so}}
\newcommand{\sym}{\mathrm{Sym}}
\newcommand{\Innl}{\operatorname{Inn}_\ell}

\title[Profinite and $p$-adic analytic Bruck loops]{Profinite and $p$-adic analytic Bruck loops as symmetric sections of congruence subgroups}

\author{Raffaello Caserta}
\address{Istituto Superiore ``Ettore Majorana'', Palermo, Italy}
\email{raffaello.caserta@majoranapa.edu.it}


\begin{document}

\begin{abstract}
We construct profinite and $p$-adic analytic Bruck loops (hence $A_\ell$-loops satisfying the left Bol and automorphic inverse properties) together with projective towers of finite Bruck loops as symmetric sections of odd-order pro-$p$ groups equipped with a continuous involution $\tau$. We first establish a global profinite Cartan--Glauberman decomposition: for every odd-order pro-$p$ group $G$ with continuous involution $\tau$, the multiplication map $L\times H\to G$ from $L=\{g\in G:\tau(g)=g^{-1}\}$ and $H=G^\tau$ is a homeomorphism, endowing $L$ with the structure of a profinite Bruck loop $a\oplus b = (ab^2a)^{1/2}$ whose quotients by open $\tau$-invariant normal subgroups form a projective system of finite Bruck loops.

For congruence subgroups $G_1\le 1+pR$ of a finite-rank $\Zp$-algebra $R$ with involution ($p\ge 3$), we parametrize $L$ globally via the Cayley transform $\kappa(x)=(1+x)(1-x)^{-1}$ ($x\in p\m$) and prove that the induced loop law $x\oplus_\kappa y$ expands as an integral formal power series in $\Z\langle\langle x,y\rangle\rangle$ governed by Catalan numbers, avoiding both square-root obstructions in $\Qp$ and Baker--Campbell--Hausdorff denominators at small primes. As our main structural result, we compute the exact congruence thresholds for the finite quotients $L_k = L/(L\cap G_k)$: the nilpotency class of the faithful left inner mapping group satisfies $\cl(\Innl(L_k))\le \lfloor(k-2)/2\rfloor$ for all $k\ge 2$, with equality whenever the reduced Lie algebra $[\bar{\m},\bar{\m}]$ is perfect and acts faithfully on $\bar{\m}$. In particular, $L_k$ is an abelian group for $k\le 3$, becomes a proper Bruck/$A_\ell$-loop at $k=4$, and attains a non-abelian $\Innl(L_k)$ precisely at $k=6$. We illustrate the theory through three classical families: the $p$-adic $2$-sphere and hyperbolic disk on $p\Zp[i]$ ($p\equiv 3\pmod 4$, of finite order $p^{2(k-1)}$), higher $p$-adic spheres $\SO_n(\Zp)_1/\SO_{n-1}(\Zp)_1$ ($n\ge 4$), and $p$-adic polar decompositions in $\GL_n(\Zp)_1/\OO_n(\Zp)_1$ and $\SL_n(\Zp)_1/\SO_n(\Zp)_1$ ($n\ge 3$).
\end{abstract}

\maketitle

\section{Introduction}

As Kikkawa first observed in~\cite{kikkawa1964local,kikkawa1975geometry}, loops arise naturally in differential geometry when one extends the affine addition of vectors to a manifold $M$ endowed with an affine connection of non-vanishing curvature. Fixing a base point $m\in M$, the geodesic composition law $a\oplus b := (\exp_a\circ P_{m,a}\circ\exp_m^{-1})(b)$ (where $P_{m,a}$ denotes parallel transport along the geodesic from $m$ to $a$) turns a normal neighborhood of $m$ into a local loop whose left inner mappings encode the holonomy of parallel transport around geodesic triangles; in this setting, curvature is directly reflected in the non-associativity of the induced loop. On a symmetric space $M=G/H$ determined by an involutive automorphism $\tau$ of a Lie group $G$, the canonical section $L = \{g\in G : \tau(g)=g^{-1}\}$ yields a (local) \emph{Bruck loop} (or \emph{$K$-loop}~\cite{kiechle2002theory,nagy2002loops}), that is, a left Bol loop satisfying the automorphic inverse property and hence in particular an \emph{$A_\ell$-loop} (left automorphic loop~\cite{bruck1956loops,kinyon2000loops,kinyon2016structure}).

In the real Archimedean setting, the loops arising in this way are infinite and often only locally defined; moreover, discretizing the real $2$-sphere in the Euclidean topology breaks the algebraic elegance of the operation~\cite{caserta2010geodesic}. Changing the underlying topology from the real to the $p$-adic one alters the algebraic structure considerably\footnote{See~\cite{caserta2011padic} for a first $1$-dimensional $p$-adic construction via M\"obius-type maps on $p\Zp$, and~\cite{caserta2017survey} for a broader geometric survey.}. Because the ultrametric norm on $\Zp$ is invariant under reduction modulo $p^k$ and open congruence subgroups of $p$-adic Lie groups are compact pro-$p$ groups~\cite{dixon1999analytic,schneider2011padic}, a global analytic section of a $p$-adic symmetric space descends, modulo powers of $p$, to a projective tower of finite $A_\ell$-loops. At the same time, from an algebraic perspective, Glauberman's classical results~\cite{glauberman1964loops,glauberman1968loops} show that finite Bruck loops of odd order are governed by groups with involutions via $a\circ b = a^{1/2}ba^{1/2}$, in sharp contrast with the $2$-power-exponent theory studied in~\cite{aschbacher2006bruck,baumeister2011bruck}.

The purpose of this paper is to develop a self-contained theory of profinite and $p$-adic analytic Bruck loops arising as global symmetric sections of odd-order pro-$p$ groups, and to determine the exact algebraic invariants of their finite congruence quotients. Our framework eliminates two obstacles that commonly affect $p$-adic geometric constructions: on the one hand, we avoid polar coordinates and quadratic square roots in $\Qp$ (which fail to exist generically in $\Qp$); on the other hand, rather than relying on the Baker--Campbell--Hausdorff formula (whose denominators $n!$ interfere with congruence filtrations at small primes $p=3,5$), we introduce a global \emph{Cayley parametrization} valid uniformly for every odd prime $p\ge 3$.

In Section~\ref{sec:pro-p-cartan} we prove a global profinite Cartan--Glauberman decomposition (Theorem~\ref{thm:pro-p-cartan}): for any pro-$p$ group $G$ ($p>2$) equipped with a continuous involution $\tau$, the set $L=\{g\in G : \tau(g)=g^{-1}\}$ is a closed section of $G/H$ ($H=G^\tau$) such that $L\times H\to G$ is a homeomorphism, $(L,\oplus)$ with $a\oplus b=(ab^2a)^{1/2}$ is a profinite Bruck (hence $A_\ell$-) loop, and its quotients modulo open $\tau$-invariant normal subgroups form a projective system of finite Bruck loops whose inverse limit recovers $L$.

In Section~\ref{sec:cayley-thresholds} we specialize to congruence subgroups $G_1\le 1+pR$ inside a finite-rank associative $\Zp$-algebra $R$ equipped with a $\Zp$-linear involution (covering both internal conjugations $\tau(g)=JgJ^{-1}$ and Cartan transpose-inverses $\tau(g)=(g^T)^{-1}$). Using the Cayley transform $\kappa(x)=(1+x)(1-x)^{-1}$, we prove that the induced loop law $x\oplus_\kappa y$ expands as an integral formal power series in $\Z\langle\langle x,y\rangle\rangle$ whose coefficients are controlled by Catalan numbers (Theorem~\ref{thm:catalan-cayley}). Computing this series explicitly up to degree~$4$ yields our main structural theorem (Theorem~\ref{thm:exact-thresholds}): for every odd prime $p\ge 3$ and every congruence level $k\ge 2$, the faithful left inner mapping group $\Innl(L_k)$ of the finite Bruck loop $L_k=L/(L\cap G_k)$ satisfies the sharp bound
\[
  \cl\bigl(\Innl(L_k)\bigr) \le \left\lfloor \frac{k-2}{2} \right\rfloor,
\]
with equality whenever the reduced Lie algebra $\bar{\h}_0 = [\bar{\m},\bar{\m}]$ is perfect and acts faithfully on $\bar{\m}$. In particular, $L_k$ is an abelian group for $k\le 3$, becomes a proper Bruck/$A_\ell$-loop precisely at $k=4$, and acquires a non-abelian left inner mapping group $\Innl(L_k)$ precisely at $k=6$ (whereas ambient commutators in $H_1/H_k$ already appear at $k=5$, lying inside the core of the action on $L_5$).

Finally, Section~\ref{sec:classical-families} applies the theory to three classical families of $p$-adic symmetric spaces: the $p$-adic $2$-sphere and hyperbolic disk on $D=p\Zp[i]$ ($p\equiv 3\pmod 4$), yielding finite Bruck/$A_\ell$-loops of exact order $p^{2(k-1)}$ with cyclic $\Innl(L_k)\cong\Z/p^{k-3}\Z$ ($k\ge 3$) and explicit cubic polynomial law modulo $p^5$ (Proposition~\ref{prop:padic-mobius}); higher $p$-adic spheres $\SO_n(\Zp)_1/\SO_{n-1}(\Zp)_1$ ($n\ge 4$), yielding finite Bruck/$A_\ell$-loops of order $p^{(n-1)(k-1)}$ attaining exact nilpotency class $\lfloor(k-2)/2\rfloor$ for every odd prime $p\ge 3$ (Theorem~\ref{thm:higher-spheres}); and $p$-adic polar decompositions in $\GL_n(\Zp)_1/\OO_n(\Zp)_1$ and $\SL_n(\Zp)_1/\SO_n(\Zp)_1$ ($n\ge 3$), likewise achieving arbitrarily high nilpotency class (Theorem~\ref{thm:polar-decomp}).

\section{Symmetric sections of odd-order \texorpdfstring{pro-$p$}{pro-p} groups}\label{sec:pro-p-cartan}

We recall the loop-theoretic notions needed below; standard references are~\cite{baer1939nets,kiechle2002theory,kinyon2000loops,nagy2002loops,sabinin1972equivalence,aschbacher2006bruck}. A \emph{loop} $(L,\oplus)$ is a set equipped with a binary operation $\oplus$ and a two-sided identity $1\in L$ such that for all $a,b\in L$ the equations $a\oplus x = b$ and $x\oplus a = b$ each have a unique solution in $L$. Denoting by $\lambda_a\colon x\mapsto a\oplus x$ the left translation by $a$, the \emph{left multiplication} (or \emph{enveloping}) group is $G(L) := \langle \lambda_a : a\in L\rangle \le \Sym(L)$, and the \emph{left inner mapping group} is the stabilizer $\Innl(L) := \Stab_{G(L)}(1) = \langle \delta_{a,b} : a,b\in L\rangle$, where $\delta_{a,b} := \lambda_{a\oplus b}^{-1}\lambda_a\lambda_b$.

A loop $(L,\oplus)$ is:
\begin{itemize}
  \item an \emph{$A_\ell$-loop} (or \emph{left automorphic loop}) if $\Innl(L)\le\Aut(L,\oplus)$;
  \item a \emph{left Bol loop} if $a\oplus(b\oplus(a\oplus c)) = (a\oplus(b\oplus a))\oplus c$ for all $a,b,c\in L$ (which implies the \emph{left inverse property}, LIP: $a^{-1}\oplus(a\oplus b)=b$);
  \item a \emph{Bruck loop} (or \emph{$K$-loop}) if it is a left Bol loop satisfying the \emph{automorphic inverse property} (AIP), namely $(a\oplus b)^{-1} = a^{-1}\oplus b^{-1}$ for all $a,b\in L$.
\end{itemize}
Every Bruck loop is in particular an $A_\ell$-loop~\cite{kiechle2002theory,nagy2002loops}. Given a group $G$, a subgroup $H\le G$, and a left transversal $L$ to $H$ in $G$ containing $1$, $L$ inherits a left-loop operation $a\oplus b := c$, where $\{c\} = abH\cap L$ (and $(L,\oplus)$ is a loop if and only if $L$ is a left transversal to all conjugates $gHg^{-1}$~\cite{baer1939nets,sabinin1972equivalence}). Writing $\delta(a,b) := (a\oplus b)^{-1}ab\in H$ for the group-theoretic left deviation in $H$, whenever $L^H=L$ every $h\in H$ acts on $(L,\oplus)$ by conjugation $a\mapsto hah^{-1}$ as an automorphism, $\delta_{a,b}$ coincides with conjugation by $\delta(a,b)$ on $L$, and the faithful quotient $D/\ker(\Ad_D|_L)$ of $D := \langle \delta(a,b) : a,b\in L\rangle\le H$ is canonically isomorphic to $\Innl(L)$.

A \emph{profinite group} is an inverse limit of finite discrete groups (equivalently, a compact, Hausdorff, totally disconnected topological group~\cite{ribes2010profinite}); it is a \emph{pro-$p$ group} if every continuous finite quotient is a finite $p$-group. Similarly, we define a \emph{profinite loop} to be a topological loop (with continuous multiplication and left/right divisions) that is isomorphic as a topological loop to an inverse limit of finite discrete loops; in particular, every profinite loop is compact, Hausdorff, and totally disconnected (cf.~\cite{herfort2011boolean} for the relation with general Boolean loops).

\begin{lemma}[Unique square roots in odd-order pro-$p$ groups]\label{lem:unique-sqrt}
Let $p$ be an odd prime and let $G$ be a pro-$p$ group. Then the squaring map $\mathrm{sq}\colon G\to G$, $x\mapsto x^2$, is a homeomorphism. Denoting its continuous inverse by $x\mapsto x^{1/2}$, one has:
\begin{enumerate}
  \item $\phi(x^{1/2}) = \phi(x)^{1/2}$ for every continuous group homomorphism $\phi\colon G\to G'$ between pro-$p$ groups;
  \item $(x^{-1})^{1/2} = (x^{1/2})^{-1}$ for all $x\in G$.
\end{enumerate}
\end{lemma}

\begin{proof}
Write $G = \varprojlim_{N\trianglelefteq_o G} G/N$ as the inverse limit of its finite $p$-group quotients. For each open normal subgroup $N\trianglelefteq_o G$, the order $|G/N| = p^{m_N}$ is odd, so $2$ is invertible modulo $p^{m_N}$; hence squaring on $G/N$ is a bijection with inverse $y\mapsto y^{(p^{m_N}+1)/2}$. Passing to the inverse limit shows that $\mathrm{sq}\colon G\to G$ is a continuous bijection from a compact space onto a Hausdorff space, therefore a homeomorphism. Properties~(1) and~(2) follow immediately from the uniqueness of square roots: $\phi(x^{1/2})^2 = \phi(x)$ and $((x^{1/2})^{-1})^2 = x^{-1}$.
\end{proof}

\begin{theorem}[Global profinite Cartan--Glauberman decomposition]\label{thm:pro-p-cartan}
Let $p$ be an odd prime, let $G$ be a pro-$p$ group, and let $\tau\in\Aut(G)$ be a continuous involutive automorphism ($\tau^2=\id_G$). Define
\[
  H := \{h\in G : \tau(h) = h\}, \qquad L := \{a\in G : \tau(a) = a^{-1}\}.
\]
Then the following hold:
\begin{enumerate}
  \item $H$ is a closed subgroup of $G$, and $L$ is a closed subset of $G$ satisfying $1\in L$, $L^{-1}=L$, $L^H=L$, and $L\cap H = \{1\}$.
  \item The multiplication map $\mu\colon L\times H\to G$, $(a,h)\mapsto ah$, is a homeomorphism, with continuous inverse
  \[
    \mu^{-1}(g) = \Bigl(\bigl(g\tau(g)^{-1}\bigr)^{1/2},\;\; \bigl(g\tau(g)^{-1}\bigr)^{-1/2}g\Bigr).
  \]
  In particular, $L$ is a continuous left transversal (closed section) to $H$ in $G$.
  \item The induced transversal operation $a\oplus b \in abH\cap L$ is given explicitly by
  \begin{equation}\label{eq:bruck-op}
    a\oplus b = (ab^2a)^{1/2} \qquad (a,b\in L),
  \end{equation}
  and makes $(L,\oplus)$ a compact, Hausdorff, totally disconnected topological Bruck loop (in fact a profinite Bruck loop by~(4), hence an $A_\ell$-loop with LIP and AIP). The squaring map $x\mapsto x^2$ is a topological loop isomorphism from $(L,\oplus)$ onto $(L,\circ)$, where $a\circ b := a^{1/2}ba^{1/2}$ is the Glauberman operation~\cite{glauberman1964loops}.
  \item For every open $\tau$-invariant normal subgroup $N\trianglelefteq_o G$, let $\pi_N\colon G\to G_N := G/N$ be the canonical projection and let $\bar{\tau}$ be the induced involution on $G_N$. Then $\pi_N(L) = \{x\in G_N : \bar{\tau}(x)=x^{-1}\}$ is isomorphic as a finite Bruck loop to the quotient $L_N := L/(L\cap N)$, and $\pi_N|_L\colon (L,\oplus)\to(L_N,\oplus)$ is a continuous loop epimorphism. Moreover, letting $\mathcal{N}_\tau$ be the cofinal filter base of open $\tau$-invariant normal subgroups of $G$, one has
  \[
    (L,\oplus) \cong \varprojlim_{N\in\mathcal{N}_\tau} (L_N,\oplus).
  \]
\end{enumerate}
\end{theorem}

\begin{proof}
(1) Since $\tau$ is continuous and $G$ is Hausdorff, $H$ and $L$ are closed. Clearly $1\in L$, $L^{-1}=L$, and for $a\in L, h\in H$ one has $\tau(hah^{-1}) = ha^{-1}h^{-1} = (hah^{-1})^{-1}$, so $L^H=L$. If $x\in L\cap H$, then $x = \tau(x) = x^{-1}$, whence $x^2 = 1$ and therefore $x=1$ by Lemma~\ref{lem:unique-sqrt}.

(2) For any $g\in G$, set $u := g\tau(g)^{-1}$. Then $\tau(u) = \tau(g)g^{-1} = u^{-1}$, so $u\in L$. Setting $a := u^{1/2}$, Lemma~\ref{lem:unique-sqrt} yields $\tau(a) = \tau(u)^{1/2} = (u^{-1})^{1/2} = a^{-1}$, hence $a\in L$. Next set $h := a^{-1}g$. We compute
\[
  \tau(h) = \tau(a)^{-1}\tau(g) = a\tau(g) = a^{-1}(a^2\tau(g)) = a^{-1}(g\tau(g)^{-1}\tau(g)) = a^{-1}g = h,
\]
so $h\in H$ and $g = ah$. To prove uniqueness, suppose $g = ah$ with $a\in L$ and $h\in H$. Then $g\tau(g)^{-1} = ah(\tau(a)\tau(h))^{-1} = ah(a^{-1}h)^{-1} = ahh^{-1}a = a^2$, forcing $a = (g\tau(g)^{-1})^{1/2}$ by Lemma~\ref{lem:unique-sqrt} and hence $h = a^{-1}g$. Because both $\mu$ and $\mu^{-1}$ are composed of continuous operations, $\mu$ is a homeomorphism.

(3) For $a,b\in L$, writing $ab = ch$ with $c\in L$ and $h\in H$, part~(2) gives
\[
  c = \bigl((ab)\tau(ab)^{-1}\bigr)^{1/2} = \bigl(ab(a^{-1}b^{-1})^{-1}\bigr)^{1/2} = (ab^2a)^{1/2},
\]
proving~\eqref{eq:bruck-op}. Clearly $1\oplus a = a\oplus 1 = a$, and $\oplus$ is continuous. To verify that $(L,\oplus)$ is a topological loop, we solve the left and right division equations explicitly. For left division, $a\oplus x = b$ holds if and only if $(ax^2a)^{1/2} = b$, equivalently $ax^2a = b^2$, i.e., $x^2 = a^{-1}b^2a^{-1}$, whose unique solution in $G$ is $x = (a^{-1}b^2a^{-1})^{1/2} = a^{-1}\oplus b$ by Lemma~\ref{lem:unique-sqrt}; since $\tau(a^{-1}b^2a^{-1}) = ab^{-2}a = (a^{-1}b^2a^{-1})^{-1}$, one has $x\in L$, and this simultaneously yields the LIP $a^{-1}\oplus(a\oplus b) = b$. For right division, $x\oplus a = b$ holds if and only if $xa^2x = b^2$; multiplying on the left and on the right by $a$ shows that this is equivalent to $(axa)^2 = ab^2a$. Since both $axa$ and $ab^2a$ lie in $L$ (as $\tau(axa) = a^{-1}x^{-1}a^{-1} = (axa)^{-1}$), uniqueness of square roots yields $axa = (ab^2a)^{1/2} = a\oplus b$, whence $x = a^{-1}(a\oplus b)a^{-1}$ is the unique solution; one checks at once that $\tau(x) = a(a\oplus b)^{-1}a = x^{-1}$, so $x\in L$, and both divisions depend continuously on $(a,b)$.

We next verify the left Bol and automorphic inverse identities directly from~\eqref{eq:bruck-op}. For all $a,b,c\in L$, one has
\[
  \bigl(a\oplus(b\oplus(a\oplus c))\bigr)^2 = a(b\oplus(a\oplus c))^2 a = ab(a\oplus c)^2 ba = aba c^2 aba.
\]
On the other hand, setting $d := a\oplus(b\oplus a)$, we obtain $d^2 = a(b\oplus a)^2 a = aba^2ba = (aba)^2$, whence $d = aba$ by Lemma~\ref{lem:unique-sqrt} (since $\tau(aba)=a^{-1}b^{-1}a^{-1}=(aba)^{-1}$ implies $aba\in L$). Therefore
\[
  \bigl((a\oplus(b\oplus a))\oplus c\bigr)^2 = d c^2 d = aba c^2 aba = \bigl(a\oplus(b\oplus(a\oplus c))\bigr)^2,
\]
and uniqueness of square roots yields $a\oplus(b\oplus(a\oplus c)) = (a\oplus(b\oplus a))\oplus c$. In addition, by~\eqref{eq:bruck-op} and Lemma~\ref{lem:unique-sqrt}(2),
\[
  (a\oplus b)^{-1} = \bigl((ab^2a)^{1/2}\bigr)^{-1} = (a^{-1}b^{-2}a^{-1})^{1/2} = a^{-1}\oplus b^{-1},
\]
establishing the AIP. Thus $(L,\oplus)$ is a profinite Bruck loop (and hence an $A_\ell$-loop, as also follows directly from $L^H=L$). Finally, for $\Phi(x) := x^2$, we have $\Phi(a\oplus b) = (a\oplus b)^2 = ab^2a = (a^2)^{1/2} b^2 (a^2)^{1/2} = \Phi(a)\circ\Phi(b)$, so $\Phi\colon(L,\oplus)\xrightarrow{\sim}(L,\circ)$ is a topological loop isomorphism.

(4) Let $N\trianglelefteq_o G$ be $\tau$-invariant. Clearly $\pi_N(L)\subseteq L(G_N) := \{x\in G_N : \bar{\tau}(x)=x^{-1}\}$. Conversely, if $gN\in L(G_N)$, then $\tau(g) \equiv g^{-1} \pmod N$, so $g\tau(g)^{-1} \equiv g^2 \pmod N$. Setting $a := (g\tau(g)^{-1})^{1/2}\in L$, we obtain $\pi_N(a)^2 = \pi_N(g)^2$ in the odd-order finite $p$-group $G_N$, whence $\pi_N(a) = \pi_N(g)$ by uniqueness of square roots. Thus $\pi_N(L) = L(G_N)$. Since $\pi_N$ is a group homomorphism, $\pi_N(a\oplus b) = \pi_N((ab^2a)^{1/2}) = (\pi_N(a)\pi_N(b)^2\pi_N(a))^{1/2} = \pi_N(a)\oplus \pi_N(b)$ by Lemma~\ref{lem:unique-sqrt}(1), so $\pi_N|_L$ is a continuous loop epimorphism with kernel $L\cap N$. Every open normal subgroup $U\trianglelefteq_o G$ contains the $\tau$-invariant open normal subgroup $U\cap\tau(U)\in\mathcal{N}_\tau$, so $\bigcap_{N\in\mathcal{N}_\tau} N = \{1\}$ and $L\cong\varprojlim_{N\in\mathcal{N}_\tau} L_N$.
\end{proof}

\section{Cayley coordinates, Catalan power series, and exact nilpotency thresholds}\label{sec:cayley-thresholds}

To obtain explicit algebraic formulas and exact congruence thresholds for finite quotients of $p$-adic Lie loops---valid uniformly for all odd primes $p\ge 3$ without interference from $n!$ denominators---we now pass from the exponential map to the Cayley transform inside an associative $\Zp$-algebra.

\subsection{The Cayley--Catalan framework}

Let $p\ge 3$ be an odd prime, and let $R$ be an associative unital $\Zp$-algebra which is free of finite rank as a $\Zp$-module (for instance $R = M_n(\Zp)$ or $M_n(\Oint)$ for a finite unramified extension $\Oint/\Zp$). Let $1 + pR$ be the first congruence subgroup of $R^\times$; it is a compact pro-$p$ group with filtration by open normal subgroups $1 + p^k R$ ($k\ge 1$).

Suppose $R$ is equipped with a $\Zp$-linear involution $\sigma\colon R\to R$ ($\sigma^2=\id_R$) of one of two types:
\begin{enumerate}
  \item[(I)] \emph{Automorphism type:} $\sigma(xy) = \sigma(x)\sigma(y)$ for all $x,y\in R$. We define the continuous group involution $\tau(g) := \sigma(g)$ on $1+pR$ and set $\theta := \sigma$ on $R$.
  \item[(II)] \emph{Anti-automorphism type:} $\sigma(xy) = \sigma(y)\sigma(x)$ for all $x,y\in R$. We define the continuous group involution $\tau(g) := \sigma(g)^{-1}$ on $1+pR$ and set $\theta := -\sigma$ on $R$.
\end{enumerate}
In both types, $\tau$ is a continuous involutive group automorphism of $1+pR$ preserving each $1+p^k R$, and $\theta\colon R\to R$ is a $\Zp$-linear involution which is an automorphism of the Lie algebra $(R,[\cdot,\cdot])$, where $[x,y] := xy - yx$. Since $2\in\Zp^\times$, $R$ decomposes into the direct sum of $\pm 1$-eigenspaces of $\theta$:
\[
  R = \h_R \oplus \m_R, \quad
  \h_R := \{x\in R : \theta(x)=x\}, \quad
  \m_R := \{x\in R : \theta(x)=-x\},
\]
satisfying $[\h_R,\h_R]\subseteq\h_R$, $[\h_R,\m_R]\subseteq\m_R$, and $[\m_R,\m_R]\subseteq\h_R$.

Define the \emph{Cayley transform} $\kappa\colon pR\to 1+pR$ and its inverse $\kappa^{-1}\colon 1+pR\to pR$ by
\begin{equation}\label{eq:cayley-def}
  \kappa(x) := (1+x)(1-x)^{-1} = 1 + 2\sum_{j=1}^\infty x^j, \qquad
  \kappa^{-1}(g) := (g-1)(g+1)^{-1}.
\end{equation}
Because $2\in\Zp^\times$ and $(1-x)^{-1} = \sum_{j\ge 0} x^j$ converges in $R$ for all $x\in pR$ (and $(g+1)^{-1} = 2^{-1}(1 + \frac{g-1}{2})^{-1}$ converges for all $g\in 1+pR$), $\kappa$ is a $p$-adic analytic diffeomorphism mapping $p^k R$ onto $1+p^k R$ for every $k\ge 1$. Moreover, a direct verification in both types~(I) and~(II) shows that
\begin{equation}\label{eq:cayley-equivariance}
  \kappa(-x) = \kappa(x)^{-1}, \qquad \tau(\kappa(x)) = \kappa(\theta(x)) \qquad (x\in pR).
\end{equation}
Consequently, $\kappa$ maps $p^k\m_R$ bijectively onto $L(1+p^k R) = \{g\in 1+p^k R : \tau(g)=g^{-1}\}$ and maps $p^k\h_R$ bijectively onto $H(1+p^k R) = (1+p^k R)^\tau$ for every $k\ge 1$.

\begin{theorem}[Integral Catalan power series in Cayley coordinates]\label{thm:catalan-cayley}
Let $p\ge 3$ be an odd prime and let $(R,\sigma,\tau,\theta)$ be as above. For all $x,y\in p\m_R$, let $a:=\kappa(x)$ and $b:=\kappa(y)$ in $L(1+pR)$. Then:
\begin{enumerate}
  \item $ab^2a = 1 + 4W(x,y)$, where
  \[
    W(x,y) := x(1-x)^{-2} + \kappa(x)\,y(1-y)^{-2}\,\kappa(x) \in pR
  \]
  expands as an integral power series in $\Z\langle\langle x,y\rangle\rangle$.
  \item The square root $(ab^2a)^{1/2} = a\oplus b$ is given by
  \[
    a\oplus b = 1 + 2S(x,y), \qquad
    S(x,y) := \sum_{m=1}^\infty (-1)^{m-1} C_{m-1}\,W(x,y)^m \in pR,
  \]
  where $C_r := \frac{1}{r+1}\binom{2r}{r}\in\Z$ denotes the $r$-th Catalan number.
  \item The induced Bruck loop operation in Cayley coordinates, $x\oplus_\kappa y := \kappa^{-1}(\kappa(x)\oplus\kappa(y))$, is given by
  \begin{equation}\label{eq:cayley-catalan-law}
    x\oplus_\kappa y = S(x,y)\bigl(1 + S(x,y)\bigr)^{-1} \in \Z\langle\langle x,y\rangle\rangle,
  \end{equation}
  a formal power series with exclusively integer coefficients, converging on $p\m_R$ for every odd prime $p\ge 3$.
  \item Modulo terms of degree~$\ge 5$ in the free associative algebra $\Z\langle\langle x,y,z,w\rangle\rangle$ (resp.\ degree~$\ge 4$ in~\eqref{eq:deg4-dev}), and hence modulo $p^5 R$ (resp.\ $p^4 R$) for $x,y,z,w\in p\m_R$, the following exact expansions hold:
  \begin{align}
    x\oplus_\kappa y &\equiv x + y - (x^2y - xyx + yx^2) \notag\\
                     &\quad + (xy^2 - 3yxy + y^2x) \pmod{p^5 R},\label{eq:deg4-op}\\
    (x\oplus_\kappa y)\oplus_\kappa z - x\oplus_\kappa(y\oplus_\kappa z) &\equiv -2\,\bigl[[x,y],\,z\bigr] \pmod{p^5 R},\label{eq:deg4-assoc}\\
    \delta\bigl(\kappa(x),\kappa(y)\bigr) &\equiv 1 + 2\,[x,y] \pmod{p^4 R},\label{eq:deg4-dev}\\
    \bigl[\delta(\kappa(x),\kappa(y)),\,\delta(\kappa(z),\kappa(w))\bigr] &\equiv 1 + 4\,\bigl[[x,y],[z,w]\bigr] \pmod{p^5 R},\label{eq:deg4-comm}
  \end{align}
  where $\delta(a,b) := (a\oplus b)^{-1}ab \in H(1+pR)$ and $[u,v]:=uv-vu$ in~\eqref{eq:deg4-assoc}--\eqref{eq:deg4-comm} denotes the algebra commutator (while on the left-hand side of~\eqref{eq:deg4-comm}, $[g_1,g_2]:=g_1g_2g_1^{-1}g_2^{-1}$ is the group commutator).
\end{enumerate}
\end{theorem}

\begin{proof}
(1) Using $\kappa(x) = 1 + 2x(1-x)^{-1}$ and $\kappa(y)^2 = (1+y)^2(1-y)^{-2} = 1 + 4y(1-y)^{-2}$, we expand
\begin{align*}
  \kappa(x)\kappa(y)^2\kappa(x)
    &= \kappa(x)^2 + 4\,\kappa(x)\,y(1-y)^{-2}\,\kappa(x) \\
    &= 1 + 4x(1-x)^{-2} + 4\,\kappa(x)\,y(1-y)^{-2}\,\kappa(x) = 1 + 4W(x,y).
\end{align*}
Since $(1-x)^{-1} = \sum_{j\ge 0} x^j\in\Z\langle\langle x\rangle\rangle$, the series $W(x,y)$ has all coefficients in $\Z$ and no constant term, so $W(x,y)\in pR$.

(2) Recall the binomial coefficient identity for $m\ge 1$:
\begin{align*}
  4^m\binom{1/2}{m} &= 4^m \frac{\tfrac{1}{2}(-\tfrac{1}{2})(-\tfrac{3}{2})\cdots(\tfrac{3}{2}-m)}{m!} \\
                    &= 2(-1)^{m-1}\frac{1}{m}\binom{2m-2}{m-1} = 2(-1)^{m-1} C_{m-1} \in 2\Z.
\end{align*}
Therefore the binomial series $(1+4W)^{1/2} = 1 + \sum_{m=1}^\infty 4^m\binom{1/2}{m} W^m = 1 + 2S(x,y)$ converges in $1+pR$ (as $W\in pR$ and $C_{m-1}\in\Z$), satisfies $(1+2S)^2 = 1+4W = ab^2a$, and hence equals $(ab^2a)^{1/2}$ by Lemma~\ref{lem:unique-sqrt}.

(3) Applying $\kappa^{-1}(g) = (g-1)(g+1)^{-1}$ to $g = 1 + 2S(x,y)$ cancels the factor of~$2$ identically:
\[
  x\oplus_\kappa y = (2S)(2+2S)^{-1} = S(1+S)^{-1} = \sum_{r=1}^\infty (-1)^{r-1} S^r \in \Z\langle\langle x,y\rangle\rangle.
\]

(4) Expanding $W(x,y)$ to degree~$4$ gives
\begin{align*}
  W(x,y) &\equiv x(1+2x+3x^2+4x^3) \\
         &\quad + (1+2x+2x^2+2x^3)\,y(1+2y+3y^2+4y^3)\,(1+2x+2x^2+2x^3) \\
         &\equiv (x+y) + 2(x^2 + xy + yx + y^2) \\
         &\quad + (3x^3 + 2x^2y + 4xyx + 4xy^2 + 2yx^2 + 4y^2x + 3y^3) \\
         &\quad + 2(2x^4 + x^3y + 2x^2yx + 2x^2y^2 + 2xyx^2 + 4xy^2x \\
         &\qquad + 3xy^3 + yx^3 + 2y^2x^2 + 3y^3x + 2y^4),
\end{align*}
modulo terms of degree~$\ge 5$. Substituting $C_0=1$, $C_1=1$, $C_2=2$, $C_3=5$ into $S = W - W^2 + 2W^3 - 5W^4 + O(5)$ and then $x\oplus_\kappa y = S - S^2 + S^3 - S^4 + O(5)$, all degree-$2$ and degree-$4$ terms cancel identically; indeed, directly in the free formal power series ring $\Z\langle\langle x,y\rangle\rangle$, the identity $\kappa(-x)=\kappa(x)^{-1}$ implies $\kappa(-x)\kappa(-y)^2\kappa(-x) = (ab^2a)^{-1}$, whence by uniqueness of formal square roots with constant term~$1$ one obtains $\kappa((-x)\oplus_\kappa(-y)) = (a\oplus b)^{-1} = \kappa(-(x\oplus_\kappa y))$, forcing $x\oplus_\kappa y$ to be an odd formal power series in $(x,y)$. This leaves precisely~\eqref{eq:deg4-op}. Applying~\eqref{eq:deg4-op} twice to compute $(x\oplus_\kappa y)\oplus_\kappa z$ and $x\oplus_\kappa(y\oplus_\kappa z)$ modulo degree~$5$ (where the degree-$4$ part again vanishes identically by parity), the cubic terms of the associator simplify to
\begin{align*}
  (x\oplus_\kappa y)\oplus_\kappa z - x\oplus_\kappa(y\oplus_\kappa z)
    &\equiv -2xyz + 2yxz + 2zxy - 2zyx \\
    &= -2\bigl((xy-yx)z - z(xy-yx)\bigr) \\
    &= -2\bigl[[x,y],z\bigr] \pmod{\deg 5},
\end{align*}
establishing~\eqref{eq:deg4-assoc}. Next, since $a\oplus b = \kappa(x\oplus_\kappa y) = 1 + 2(x+y) + 2(x+y)^2 + O(3)$ and $ab = (1+2x+2x^2+O(3))(1+2y+2y^2+O(3)) = 1 + 2(x+y) + 2x^2 + 4xy + 2y^2 + O(3)$, we obtain
\[
  \delta(a,b) = (a\oplus b)^{-1}ab \equiv 1 + (4xy - 2xy - 2yx) = 1 + 2[x,y] \pmod{\deg 3}.
\]
Directly in $\Z\langle\langle x,y\rangle\rangle$, substituting $(x,y)\mapsto(-x,-y)$ yields
\begin{align*}
  \delta(a^{-1},b^{-1})
    &= (a^{-1}\oplus b^{-1})^{-1}a^{-1}b^{-1} = (a\oplus b)a^{-1}b^{-1} \\
    &= (a\oplus b)^{-1}(ab^2a)a^{-1}b^{-1} = (a\oplus b)^{-1}ab = \delta(a,b),
\end{align*}
so $\delta(\kappa(x),\kappa(y))$ is an even formal power series in $(x,y)$; consequently its degree-$3$ term vanishes identically and~\eqref{eq:deg4-dev} holds modulo degree~$4$ (i.e., modulo $p^4 R$). Finally, for any $u=1+2A+O(4)$ and $v=1+2B+O(4)$ with $A=[x,y]\in p^2 R$ and $B=[z,w]\in p^2 R$, expanding $uvu^{-1}v^{-1}$ modulo $p^5 R$ immediately gives $uvu^{-1}v^{-1} \equiv 1 + 4[A,B] \pmod{p^5 R}$, proving~\eqref{eq:deg4-comm}.
\end{proof}

\subsection{Exact congruence thresholds and nilpotency class of \texorpdfstring{$\Innl(L_k)$}{Inn\_l(L\_k)}}

We now consider a closed $\tau$-invariant subgroup $G_1\le 1+pR$ (such as $1+pR$ itself, or a classical congruence subgroup like $\SO_n(\Zp)_1$) whose Cayley logarithm $\mathcal{L} := p^{-1}\kappa^{-1}(G_1)\subseteq R$ is a $\Zp$-direct summand of $R$ satisfying $\theta(\mathcal{L}) = \mathcal{L}$ and $\kappa(p^k\mathcal{L}) = G_k := G_1\cap(1+p^k R)$ for every $k\ge 1$. Write the induced eigenspace decomposition as
\[
  \mathcal{L} = \h \oplus \m, \qquad \h := \mathcal{L}\cap\h_R, \qquad \m := \mathcal{L}\cap\m_R,
\]
and set $H_k := G_k^\tau = \kappa(p^k\h)$ and $L := \{a\in G_1 : \tau(a)=a^{-1}\} = \kappa(p\m)$. Notice that $\mathcal{L}$ is automatically a $\Zp$-Lie subalgebra of $(R,[\cdot,\cdot])$ with $[\m,\m]\subseteq\h$, $[\h,\m]\subseteq\m$, and $[\h,\h]\subseteq\h$: indeed, for any $X,Y\in\mathcal{L}$ and every $m\ge 1$, the group commutator $[\kappa(p^m X),\kappa(p^m Y)]\equiv 1+4p^{2m}[X,Y]\pmod{p^{3m}R}$ lies in $G_{2m} = \kappa(p^{2m}\mathcal{L})$; writing it as $\kappa(p^{2m}Z_m) = 1+2p^{2m}Z_m+O(p^{4m}R)$ with $Z_m\in\mathcal{L}$ gives $2[X,Y]\equiv Z_m\pmod{p^m R}$ for all $m\ge 1$, whence $[X,Y]\in \bigcap_{m\ge 1}(\mathcal{L}+p^m R) = \mathcal{L}$ since $2\in\Zp^\times$ and $\mathcal{L}$ is closed in $R$.

For each $k\ge 2$, let $L_k := L/(L\cap G_k)$ be the finite Bruck loop of Theorem~\ref{thm:pro-p-cartan}(4). From the non-commutative identity $\kappa(x)-\kappa(y) = 2(1-x)^{-1}(x-y)(1-y)^{-1}$ (for $x,y\in pR$, where $2\in\Zp^\times$ and $(1-x),(1-y)\in 1+pR$), one has $\kappa(x)\equiv\kappa(y)\pmod{p^k R}$ if and only if $x\equiv y\pmod{p^k R}$. Under the Cayley diffeomorphism $\kappa\colon p\m\xrightarrow{\sim} L$, $L_k$ is therefore parametrized bijectively by $p\m/p^k\m \cong (\Z/p^{k-1}\Z)^d$, where $d := \rk_{\Zp}(\m)$, so
\begin{equation}\label{eq:finite-order}
  |L_k| = p^{d(k-1)}.
\end{equation}
Denote by $\bar{\m} := \m/p\m$ and $\bar{\h} := \h/p\h$ the reductions modulo $p$, and by $\bar{\h}_0 := [\bar{\m},\bar{\m}]\subseteq\bar{\h}$ the $\Fp$-subspace spanned by all commutators $[X,Y]$ ($X,Y\in\bar{\m}$). Because $[\bar{\h},\bar{\m}]\subseteq\bar{\m}$ and the Jacobi identity gives $[[\bar{\m},\bar{\m}],[\bar{\m},\bar{\m}]]\subseteq[[\bar{\m},\bar{\m}],\bar{\h}]\subseteq[\bar{\m},\bar{\m}]$, $\bar{\h}_0$ is a Lie subalgebra of $\bar{\h}$ acting on $\bar{\m}$ by $\mathrm{ad}(A)(Z) = [A,Z]$.

Recall that for any group $D$, its lower central series is defined by $\gamma_1(D) := D$ and $\gamma_{c+1}(D) := [\gamma_c(D), D]$, and $D$ has nilpotency class $\cl(D) \le c$ if and only if $\gamma_{c+1}(D) = \{1\}$ (with $\cl(\{1\}) = 0$, abelian nontrivial groups having class~$1$, and non-abelian nilpotent groups having class~$\ge 2$).

\begin{theorem}[Exact congruence thresholds and nilpotency class]\label{thm:exact-thresholds}
Let $p\ge 3$ be an odd prime, let $G_1\le 1+pR$ be a $\tau$-invariant closed subgroup as above with $\mathcal{L}=\h\oplus\m$, and let $L_k = L/(L\cap G_k)$ ($k\ge 2$) be the finite Bruck/$A_\ell$-loop of order $p^{(k-1)\rk_{\Zp}(\m)}$. Denote by $D_1 := \langle \delta(a,b) : a,b\in L\rangle \le H_1$ the ambient deviation subgroup, and by $\Innl(L_k) \cong D_1 / \ker(\Ad|_{L_k})$ the faithful left inner mapping group of $L_k$. Then:
\begin{enumerate}
  \item For every $c\ge 1$, $\gamma_c(D_1)\le H_{2c}$, and every element of $H_{2c}$ acts trivially by conjugation on $L_{2c+1}$. Consequently, for every $k\ge 2$,
  \begin{equation}\label{eq:nilp-upper}
    \cl\bigl(\Innl(L_k)\bigr) \le \left\lfloor \frac{k-2}{2} \right\rfloor.
  \end{equation}
  \item \emph{(Non-associativity threshold)} For $k\in\{2,3\}$, $L_k$ is an abelian group ($\Innl(L_k)=\{1\}$), even though $D_1/(D_1\cap G_3)$ is nontrivial whenever $\bar{\h}_0\ne 0$. At the minimal threshold $k=4$, $L_4$ is a proper finite Bruck/$A_\ell$-loop if and only if $[\bar{\h}_0,\bar{\m}] = [[\bar{\m},\bar{\m}],\bar{\m}] \ne 0$, in which case $L_k$ is proper for every $k\ge 4$.
  \item \emph{(Non-abelian $\Innl$ threshold)} For $k\le 5$, $\Innl(L_k)$ is abelian, even though $\gamma_2(D_1)/(\gamma_2(D_1)\cap G_5)$ is nontrivial whenever $[\bar{\h}_0,\bar{\h}_0]\ne 0$. At the minimal threshold $k=6$, $\Innl(L_6)$ is non-abelian if and only if $[[\bar{\h}_0,\bar{\h}_0],\bar{\m}]\ne 0$, in which case $\Innl(L_k)$ is non-abelian for every $k\ge 6$.
  \item \emph{(Exact equality for all classes)} Suppose that the reduced Lie algebra $\bar{\h}_0 = [\bar{\m},\bar{\m}]$ is \emph{perfect} ($[\bar{\h}_0,\bar{\h}_0] = \bar{\h}_0 \ne 0$) and acts faithfully on $\bar{\m}$ (i.e., $\{A\in\bar{\h}_0 : [A,\bar{\m}]=0\} = \{0\}$). Then for every $k\ge 2$,
  \begin{equation}\label{eq:nilp-exact}
    \cl\bigl(\Innl(L_k)\bigr) = \left\lfloor \frac{k-2}{2} \right\rfloor,
  \end{equation}
  the left inner mapping group $\Innl(L)$ of $L=\varprojlim_k L_k$ is non-solvable, and its closure $\overline{\Innl(L)}$ inside the compact group $H_1/\ker(\Ad_{H_1}|_L)$ (equivalently, inside $\Aut_{\mathrm{top}}(L) \cong \varprojlim_k \Aut(L_k)$) is canonically isomorphic to $\varprojlim_k \Innl(L_k)$ and is a non-solvable $p$-adic analytic pro-$p$ group.
\end{enumerate}
\end{theorem}

\begin{proof}
(1) By Theorem~\ref{thm:catalan-cayley}(4), every generator $\delta(\kappa(x),\kappa(y))$ ($x,y\in p\m$) satisfies $\delta(\kappa(x),\kappa(y)) \equiv 1 + 2[x,y] \pmod{p^4 R}$, hence lies in $H_2 = H_1\cap (1+p^2 R)$. Therefore $D_1\le H_2$. Since $[1+p^r R, 1+p^s R]\le 1+p^{r+s} R$ for all $r,s\ge 1$, induction on $c\ge 1$ gives $\gamma_c(D_1)\le H_{2c}$.

Now take any $h = 1 + p^r A \in H_r$ ($r\ge 1$, $A\in R$) and any $a = \kappa(x) = 1 + 2x + O(p^2 R)\in L$ ($x\in p\m$). Since $a-1 = 2x + O(p^2 R)$ and $h^{-1} \equiv 1 \pmod{p^r R}$, we compute
\begin{equation}\label{eq:conj-action}
  hah^{-1} - a = (ha - ah)h^{-1} = p^r [A, a-1] h^{-1} \equiv 2p^r [A, x] \pmod{p^{r+2}R},
\end{equation}
valid for every $r\ge 1$. Because $x\in p\m$, the leading commutator term $2p^r[A,x]$ lies in $p^{r+1}R$, so $hah^{-1}\equiv a\pmod{p^{r+1}R}$ for all $a\in L$. In particular, every subgroup $H_r$ acts trivially on $L_{r+1}$; taking $r = 2c$ shows that $H_{2c}$ (and hence $\gamma_c(D_1)$) acts trivially on $L_{2c+1}$. Setting $c := \lfloor(k-2)/2\rfloor + 1$, one has $2c+1 \ge k$, whence $\gamma_c(\Innl(L_k)) = \{1\}$ and therefore $\cl(\Innl(L_k)) \le c-1 = \lfloor(k-2)/2\rfloor$.

(2) By~\eqref{eq:nilp-upper}, for $k\in\{2,3\}$ we have $\cl(\Innl(L_k)) = 0$, i.e., $\Innl(L_k)=\{1\}$, so $L_k$ is a group; moreover~\eqref{eq:deg4-op} reduces modulo $p^3 R$ to $x\oplus_\kappa y \equiv x+y\pmod{p^3 R}$, proving that $L_2$ and $L_3$ are abelian groups. Whenever $\bar{\h}_0\ne 0$, choosing $X,Y\in\m$ with $[X,Y]\not\equiv 0\pmod p$ yields $\delta(\kappa(pX),\kappa(pY))\equiv 1+2p^2[X,Y]\not\equiv 1\pmod{p^3 R}$ by~\eqref{eq:deg4-dev}, so $D_1/(D_1\cap G_3)$ is nontrivial even though it acts trivially on $L_3$. At $k=4$, for $x=pX$, $y=pY$, $z=pZ$ ($X,Y,Z\in\m$), equation~\eqref{eq:deg4-assoc} gives
\[
  (x\oplus_\kappa y)\oplus_\kappa z - x\oplus_\kappa(y\oplus_\kappa z) \equiv -2p^3\bigl[[X,Y],Z\bigr] \pmod{p^4 R}.
\]
Since $-2\in\Zp^\times$, $L_4$ is non-associative if and only if $[[\bar{\m},\bar{\m}],\bar{\m}] = [\bar{\h}_0,\bar{\m}] \ne 0$. When this holds, the loop epimorphism $L_k\twoheadrightarrow L_4$ implies that $L_k$ is non-associative for every $k\ge 4$.

(3) By~\eqref{eq:nilp-upper}, for $k\le 5$ one has $\cl(\Innl(L_k))\le 1$, so $\Innl(L_k)$ is abelian. Whenever $[\bar{\h}_0,\bar{\h}_0]\ne 0$, choosing $X,Y,Z,W\in\m$ with $[[X,Y],[Z,W]]\not\equiv 0\pmod p$ yields $[\delta(\kappa(pX),\kappa(pY)),\delta(\kappa(pZ),\kappa(pW))]\equiv 1+4p^4[[X,Y],[Z,W]]\not\equiv 1\pmod{p^5 R}$ by~\eqref{eq:deg4-comm}, so $\gamma_2(D_1)/(\gamma_2(D_1)\cap G_5)$ is nontrivial while lying in the kernel of the action on $L_5$. At $k=6$, observe that by~\eqref{eq:deg4-dev} every element of $D_1$ takes the form $1 + p^2 A + O(p^3 R)$ with $A\bmod p\in\bar{\h}_0$, whence every element of $\gamma_2(D_1) = [D_1,D_1]$ takes the form $h = 1 + p^4 C + O(p^5 R)$ with $C\bmod p\in [\bar{\h}_0,\bar{\h}_0]$; conversely, every elementary commutator $4[A,B]$ ($A=[X,Y]$, $B=[Z,W]$) is realized modulo $p^5 R$ via~\eqref{eq:deg4-comm}. Applying~\eqref{eq:conj-action} with $r=5$ shows that $H_5$ acts trivially on $L_6$, so the action of $\gamma_2(D_1)$ on $L_6$ factors exclusively through its image $[\bar{\h}_0,\bar{\h}_0]\subseteq H_4/H_5$. Applying~\eqref{eq:conj-action} with $r=4$ to $h = 1 + p^4 C + O(p^5 R)\in\gamma_2(D_1)$ and $a = \kappa(pU)$ ($U\in\m$) gives
\[
  h a h^{-1} \equiv a + 2p^5 [C, U] \pmod{p^6 R}.
\]
Since $2\in\Zp^\times$, $\gamma_2(\Innl(L_6))$ is nontrivial if and only if $[[\bar{\h}_0,\bar{\h}_0],\bar{\m}] \ne 0$. Whenever this holds, the surjection $\Innl(L_k)\twoheadrightarrow\Innl(L_6)$ ensures that $\Innl(L_k)$ is non-abelian for every $k\ge 6$.

(4) Assume that $\bar{\h}_0 = [\bar{\m},\bar{\m}]$ is perfect and acts faithfully on $\bar{\m}$. We prove by induction on $c\ge 1$ that the image of $\gamma_c(D_1)$ under the homomorphism $H_{2c}\to p^{2c}\h/p^{2c+1}\h \cong \bar{\h}$ (sending $1+p^{2c}A + O(p^{2c+1}R)\mapsto A\bmod p$) contains the entire subspace $\bar{\h}_0$. For $c=1$,~\eqref{eq:deg4-dev} shows that $\delta(\kappa(pX),\kappa(pY))\equiv 1 + 2p^2[X,Y]\pmod{p^3 R}$ maps to $2[X,Y]\bmod p$; since $2\in\Fp^\times$ and the elements $[X,Y]$ span $\bar{\h}_0$, the image of $D_1$ in $\bar{\h}$ contains $\bar{\h}_0$. Assuming the claim holds for $c\ge 1$, take any $A,B\in\bar{\h}_0$ and choose $g_1 = 1 + p^{2c}\tilde{A} + O(p^{2c+1}R)\in\gamma_c(D_1)$ and $g_2 = 1 + p^2\tilde{B} + O(p^3 R)\in D_1$ with $\tilde{A}\equiv A$, $\tilde{B}\equiv B\pmod p$. Their group commutator satisfies
\[
  [g_1,g_2] = g_1g_2g_1^{-1}g_2^{-1} \equiv 1 + p^{2c+2}[\tilde{A},\tilde{B}] \pmod{p^{2c+3}R},
\]
which maps to $[A,B]\in\bar{\h}_0$. Because $[\bar{\h}_0,\bar{\h}_0] = \bar{\h}_0$, the image of $\gamma_{c+1}(D_1)$ in $\bar{\h}$ contains $\bar{\h}_0$, completing the induction.

Now let $k\ge 2$ and set $c := \lfloor(k-2)/2\rfloor$. If $c=0$,~\eqref{eq:nilp-exact} holds trivially. If $c\ge 1$, then $k\ge 2c+2$. Choose any nonzero $A\in\bar{\h}_0$. By faithfulness of the action of $\bar{\h}_0$ on $\bar{\m}$, there exists $Z\in\m$ such that $[A,Z]\not\equiv 0\pmod p$. Picking $h = 1 + p^{2c}\tilde{A} + O(p^{2c+1}R)\in\gamma_c(D_1)$ with $\tilde{A}\equiv A\pmod p$, equation~\eqref{eq:conj-action} (with $r=2c$) yields
\[
  h\kappa(pZ)h^{-1} \equiv \kappa(pZ) + 2p^{2c+1}[\tilde{A},Z] \not\equiv \kappa(pZ) \pmod{p^{2c+2}R}.
\]
Hence $h$ acts nontrivially on $L_{2c+2}$, and a fortiori on $L_k$, so $\gamma_c(\Innl(L_k))\ne\{1\}$. Together with~\eqref{eq:nilp-upper}, this proves $\cl(\Innl(L_k)) = c = \lfloor(k-2)/2\rfloor$.

Finally, in the inverse limit $L = \varprojlim_k L_k$, the same induction applies to the derived series of $D_1$: the $r$-th derived subgroup $D_1^{(r)}$ lies in $H_{2^{r+1}}$ and maps onto $\bar{\h}_0$ modulo $p^{2^{r+1}+1}$, hence acts nontrivially on $L_{2^{r+1}+2}$ for every $r\ge 0$. Thus $\Innl(L)$ is non-solvable. Setting $K_k := \ker(\Ad_{H_1}|_{L_k})$, each $K_k$ is an open normal subgroup of $H_1$ (as $H_k\le K_k$) with $\bigcap_k K_k = \ker(\Ad_{H_1}|_L)$ and $D_1 K_k = \overline{D_1}K_k$, whence the closure $\overline{\Innl(L)}$ of $\Innl(L)$ in $H_1/\ker(\Ad_{H_1}|_L)$ satisfies $\overline{\Innl(L)} \cong \overline{D_1}/\ker(\Ad_{\overline{D_1}}|_L) \cong \varprojlim_k \overline{D_1}/(\overline{D_1}\cap K_k) = \varprojlim_k \Innl(L_k)$. Since $R$ is free of finite rank $m$ over $\Zp$, left regular multiplication embeds $1+pR$ as a closed subgroup of $\GL_m(\Zp)$; therefore $H_1$, its quotient $H_1/\ker(\Ad_{H_1}|_L)$, and the closed subgroup $\overline{\Innl(L)}$ are all $p$-adic analytic pro-$p$ groups by~\cite[Theorem~9.6]{dixon1999analytic}.
\end{proof}

\begin{remark}
It is worth emphasizing that the conditions $[\bar{\h}_0,\bar{\m}]\ne 0$ and $[[\bar{\h}_0,\bar{\h}_0],\bar{\m}]\ne 0$ in Theorem~\ref{thm:exact-thresholds}(2)--(3) are sharp at the minimal congruence levels $k=4$ and $k=6$: if $[[\m,\m],\m]$ is nonzero in $\m$ but divisible by $p$, then $L_4$ remains associative while non-associativity is merely delayed to a higher congruence level $k\ge 5$.
\end{remark}

\section{Explicit classical families of finite and \texorpdfstring{$p$-adic}{p-adic} Bruck loops}\label{sec:classical-families}

We now apply Theorems~\ref{thm:pro-p-cartan}, \ref{thm:catalan-cayley}, and~\ref{thm:exact-thresholds} to three classical families of $p$-adic symmetric spaces, obtaining explicit infinite towers of finite Bruck/$A_\ell$-loops and clarifying the arithmetic and topological subtleties of the rank-$1$ construction.

\subsection{The \texorpdfstring{$p$-adic $2$-sphere}{p-adic 2-sphere} and hyperbolic disk (rank~\texorpdfstring{$1$}{1}, abelian \texorpdfstring{$\Innl$}{Inn\_l})}

Let $p$ be a prime with $p\equiv 3\pmod 4$, so that $-1$ is not a square in $\Fp$ and $K := \Qp(i)$ ($i^2=-1$) is the unique unramified quadratic extension of $\Qp$, with ring of integers $\Oint = \Zp[i]$, maximal ideal $D := p\Oint = \{\xi\in K : |\xi|_p < 1\}$, residue field $\Oint/p\Oint\cong\F_{p^2}$, and Galois conjugation $\overline{a+ib} = a-ib$~\cite{robert2000course,schikhof1984ultrametric}.

Consider the Pauli matrix embedding of $\Qp^3$ into $M_2(K)$:
\[
  [a,b,c] := \begin{pmatrix} c & a-ib \\ a+ib & -c \end{pmatrix}, \qquad [a,b,c]^2 = (a^2+b^2+c^2)I.
\]
Recall that over $\Qp$ ($p>2$), the standard ternary quadratic form $q_0(a,b,c)=a^2+b^2+c^2$ is \emph{isotropic} (its Hilbert symbol satisfies $(-1,-1)_p = +1$~\cite{kitaoka1993arithmetic}), so $\SO_3(\Qp)\cong\PGL_2(\Qp)$ is a totally disconnected $p$-adic Lie group (whose connected components are singletons); its derived subgroup $\SO_3(\Qp)'\cong\PSL_2(\Qp)$ has index~$4$ as the kernel of the spinor norm $\SO_3(\Qp)\to\Qp^\times/(\Qp^\times)^2$. By contrast, for $p\equiv 3\pmod 4$ the anisotropic ternary form $q_{\mathrm{an}}(a,b,c) = a^2+b^2+pc^2$ has compact special orthogonal group $\SO(q_{\mathrm{an}})\cong\mathbb{H}_p^\times/\Qp^\times$, where $\mathbb{H}_p$ denotes the unique quaternion division algebra over $\Qp$. Working inside $R := M_2(\Oint)$, the first congruence unitary subgroups provide compact pro-$p$ groups directly.

Fix a sign $\varepsilon\in\{+1,-1\}$, and consider the $\Zp$-algebra $R = M_2(\Oint)$ with the conjugation involution $\sigma(M) := \sigma_z M \sigma_z^{-1}$ of type~(I), where $\sigma_z = \mathrm{diag}(1,-1)$. Let $J_\varepsilon := \mathrm{diag}(1,\varepsilon)$ and let $G_1^{(\varepsilon)}$ be the congruence subgroup
\[
  G_1^{(\varepsilon)} := \left\{ M = \begin{pmatrix} \alpha & \beta \\ -\varepsilon\overline{\beta} & \overline{\alpha} \end{pmatrix} \in I + pM_2(\Oint) : \det M = \alpha\overline{\alpha} + \varepsilon\beta\overline{\beta} = 1 \right\},
\]
equivalently $\{M\in I+pM_2(\Oint) : M^* J_\varepsilon M = J_\varepsilon,\; \det M = 1\}$ (where $M^* := \overline{M}^T$). Then $G_1^{(+1)} = \SU_2(\Oint)_1$ (isomorphic to $\SO_3(\Zp)_1$ via the adjoint action on Pauli matrices, since $-I\notin\SU_2(\Oint)_1$) and $G_1^{(-1)} = \SU(1,1;\Oint)_1$ (the $p$-adic hyperbolic counterpart). Under the Cayley transform $\kappa(X) = (I+X)(I-X)^{-1}$, the corresponding Lie algebra $\mathcal{L}^{(\varepsilon)} = \h\oplus\m^{(\varepsilon)}$ consists of
\[
  \h = \left\{ \begin{pmatrix} it & 0 \\ 0 & -it \end{pmatrix} : t\in\Zp \right\}, \qquad
  \m^{(\varepsilon)} = \left\{ \begin{pmatrix} 0 & \eta \\ -\varepsilon\overline{\eta} & 0 \end{pmatrix} : \eta\in\Oint \right\}.
\]
Indeed, for $X\in pM_2(\Oint)$, writing $M=\kappa(X)=(I+X)(I-X)^{-1}$, the unitary condition $M^* J_\varepsilon M = J_\varepsilon$ is equivalent to $(I+X^*)J_\varepsilon(I+X) = (I-X^*)J_\varepsilon(I-X)$, i.e., $X^* J_\varepsilon + J_\varepsilon X = 0$ (since $2\in\Zp^\times$), while $\det(I+X)=\det(I-X)$ is equivalent to $\tr(X)=0$. Hence $\kappa(p^k\mathcal{L}^{(\varepsilon)}) = G_k^{(\varepsilon)}$ for every $k\ge 1$, avoiding polar coordinates and square roots $\sqrt{\eta\overline{\eta}}\in\Qp$ (which fail to exist generically, e.g., for $p=3$ and $\eta = 3(1+i)$, where $\eta\overline{\eta}=18$ is not a square in $\Q_3$). Explicitly, for $X = \bigl(\begin{smallmatrix} 0 & \eta \\ -\varepsilon\overline{\eta} & 0\end{smallmatrix}\bigr)\in p\m^{(\varepsilon)}$ ($\eta\in D = p\Oint$), one has
\[
  \kappa(X) = \frac{1-\varepsilon\eta\overline{\eta}}{1+\varepsilon\eta\overline{\eta}}
  \begin{pmatrix} 1 & \xi \\ -\varepsilon\overline{\xi} & 1 \end{pmatrix},
  \qquad \xi(\eta) := \frac{2\eta}{1-\varepsilon\eta\overline{\eta}} \in D.
\]
The coordinate change $\eta\mapsto\xi(\eta)$ is a $p$-adic analytic bijection $D\xrightarrow{\sim} D$ with explicit inverse $\eta(\xi) = \xi\bigl(1+(1+\varepsilon\xi\overline{\xi})^{1/2}\bigr)^{-1}$ (since $1+\varepsilon\xi\overline{\xi}\in 1+p^2\Zp$), and is in fact a $p$-adic isometry: for all $\eta_1,\eta_2\in D$,
\[
  \xi(\eta_1)-\xi(\eta_2) = \frac{2\bigl((\eta_1-\eta_2)+\varepsilon\eta_1\eta_2(\overline{\eta}_1-\overline{\eta}_2)\bigr)}{(1-\varepsilon\eta_1\overline{\eta}_1)(1-\varepsilon\eta_2\overline{\eta}_2)},
\]
whence $|\xi(\eta_1)-\xi(\eta_2)|_p = |\eta_1-\eta_2|_p$ by the ultrametric inequality, so $\eta\mapsto\xi(\eta)$ induces a bijection modulo $p^k\Oint$ for every $k\ge 1$. Moreover, $G_1^{(\varepsilon)}$ intersects the scalar matrices trivially (if $\lambda I\in G_1^{(\varepsilon)}$, then $\lambda^2=1$ and $\lambda\equiv 1\pmod p$, forcing $\lambda=1$), so $G_1^{(\varepsilon)}$ embeds into $\PGL_2(\Oint)$, where $L^{(\varepsilon)} = \kappa(p\m^{(\varepsilon)})$ is represented projectively by the matrices
\[
  M_\varepsilon(\xi) := \begin{bmatrix} 1 & \xi \\ -\varepsilon\overline{\xi} & 1 \end{bmatrix} \in \PGL_2(\Oint) \qquad (\xi\in D = p\Oint).
\]
Multiplying two such matrices gives
\begin{align*}
  M_\varepsilon(\xi_1)M_\varepsilon(\xi_2)
  &= \begin{bmatrix} 1-\varepsilon\xi_1\overline{\xi}_2 & \xi_1+\xi_2 \\ -\varepsilon(\overline{\xi}_1+\overline{\xi}_2) & 1-\varepsilon\overline{\xi}_1\xi_2 \end{bmatrix} \\
  &= M_\varepsilon(\xi_1\oplus_\varepsilon \xi_2)
    \begin{bmatrix} 1-\varepsilon\xi_1\overline{\xi}_2 & 0 \\ 0 & 1-\varepsilon\overline{\xi}_1\xi_2 \end{bmatrix},
\end{align*}
where the induced loop operation on $D = p\Oint$ is the $p$-adic M\"obius addition~\cite{ungar1994holomorphic,caserta2017survey}
\begin{equation}\label{eq:padic-mobius}
  \xi_1\oplus_\varepsilon \xi_2 := \frac{\xi_1+\xi_2}{1-\varepsilon\overline{\xi}_1\xi_2} \qquad (\xi_1,\xi_2\in D),
\end{equation}
and the left deviation acts on $D$ by multiplication by the unitary factor
\begin{equation}\label{eq:padic-mobius-dev}
  \delta_{\xi_1,\xi_2}(\xi) = \frac{1-\varepsilon\xi_1\overline{\xi}_2}{1-\varepsilon\overline{\xi}_1\xi_2}\,\xi.
\end{equation}
Unlike the Archimedean complex setting---where the spherical law $\varepsilon=+1$ develops singularities on the unit circle whenever $\overline{\xi}_1\xi_2 = 1$, while Ungar's hyperbolic addition $\varepsilon=-1$ preserves the open unit disk $|\xi|<1$---in the ultrametric norm one has $|1-\varepsilon\overline{\xi}_1\xi_2|_p = 1$ for both signs $\varepsilon=\pm 1$, so both laws define global analytic Bruck loops on $D = p\Oint$.

\begin{proposition}[Finite $p$-adic M\"obius Bruck loops]\label{prop:padic-mobius}
Let $p\equiv 3\pmod 4$, $\varepsilon\in\{+1,-1\}$, and $k\ge 2$. Let $L_k^{(\varepsilon)} := p\Oint/p^k\Oint$ equipped with the reduction modulo $p^k$ of~\eqref{eq:padic-mobius}. Then:
\begin{enumerate}
  \item $L_k^{(\varepsilon)}$ is a finite Bruck/$A_\ell$-loop of order $p^{2(k-1)}$ (in particular, $|L_4^{(\varepsilon)}| = p^6$).
  \item Modulo $p^5$, the operation takes the explicit cubic polynomial form
  \[
    \xi_1\oplus_\varepsilon\xi_2 \equiv (\xi_1+\xi_2)\bigl(1 + \varepsilon\overline{\xi}_1\xi_2\bigr) \pmod{p^5\Oint}.
  \]
  \item For $k\in\{2,3\}$, $L_k^{(\varepsilon)}$ is an abelian group isomorphic to $(\Z/p^{k-1}\Z)^2$. For every $k\ge 4$, $L_k^{(\varepsilon)}$ is a proper non-associative, non-commutative Bruck/$A_\ell$-loop (without RIP), whose faithful left inner mapping group $\Innl(L_k^{(\varepsilon)})$ is cyclic of exact order $p^{k-3}$.
\end{enumerate}
\end{proposition}

\begin{proof}
(1) Since $\Oint/p\Oint\cong\F_{p^2}$ has $p^2$ elements, $\rk_{\Zp}(\m^{(\varepsilon)}) = 2$ and $|p\Oint/p^k\Oint| = (p^2)^{k-1} = p^{2(k-1)}$ by~\eqref{eq:finite-order}. By Theorem~\ref{thm:pro-p-cartan}(4), $L_k^{(\varepsilon)}$ is a finite Bruck/$A_\ell$-loop.

(2) For $\xi_1,\xi_2\in p\Oint$, one has $\overline{\xi}_1\xi_2\in p^2\Oint$, whence $(\overline{\xi}_1\xi_2)^2\in p^4\Oint$ and $(\xi_1+\xi_2)(\overline{\xi}_1\xi_2)^n\in p^{2n+1}\Oint\subseteq p^5\Oint$ for all $n\ge 2$. Expanding the geometric series $(1-\varepsilon\overline{\xi}_1\xi_2)^{-1} = 1 + \varepsilon\overline{\xi}_1\xi_2 + \sum_{n\ge 2}(\varepsilon\overline{\xi}_1\xi_2)^n$ gives $\xi_1\oplus_\varepsilon\xi_2 \equiv (\xi_1+\xi_2)(1+\varepsilon\overline{\xi}_1\xi_2)\pmod{p^5\Oint}$.

(3) One checks directly in $\m^{(\varepsilon)}$ that for $X_1 = \bigl(\begin{smallmatrix} 0 & 1 \\ -\varepsilon & 0\end{smallmatrix}\bigr)$ and $X_2 = \bigl(\begin{smallmatrix} 0 & i \\ \varepsilon i & 0\end{smallmatrix}\bigr)$,
\begin{align*}
  [X_1,X_2] &= \begin{pmatrix} 2\varepsilon i & 0 \\ 0 & -2\varepsilon i \end{pmatrix} \in \h, \\
  \bigl[[X_1,X_2], X_1\bigr] &= \begin{pmatrix} 0 & 4\varepsilon i \\ 4i & 0 \end{pmatrix} = 4\varepsilon X_2 \not\equiv 0 \pmod p.
\end{align*}
Hence $[\bar{\h}_0,\bar{\m}^{(\varepsilon)}]\ne 0$ while $\bar{\h}_0 = \Fp\cdot\mathrm{diag}(i,-i)$ is $1$-dimensional (abelian). By Theorem~\ref{thm:exact-thresholds}(2), $L_k^{(\varepsilon)}$ is an abelian group for $k\le 3$ and a proper non-associative Bruck/$A_\ell$-loop for all $k\ge 4$. Explicitly, for $\xi_1=p$, $\xi_2=pi$, and $\xi_3=p$, one obtains $(\xi_1\oplus_\varepsilon\xi_2)\oplus_\varepsilon\xi_3 - \xi_1\oplus_\varepsilon(\xi_2\oplus_\varepsilon\xi_3) \equiv -2\varepsilon p^3 i \not\equiv 0\pmod{p^4\Oint}$, while $\xi_1\oplus_\varepsilon\xi_2 - \xi_2\oplus_\varepsilon\xi_1 \equiv \varepsilon(\xi_1+\xi_2)(\overline{\xi}_1\xi_2-\xi_1\overline{\xi}_2)\not\equiv 0\pmod{p^4\Oint}$ (non-commutativity) and $(\xi_1\oplus_\varepsilon\xi_2)\oplus_\varepsilon(-\xi_2) \equiv p - 2\varepsilon p^3 \not\equiv p\pmod{p^4\Oint}$ (failure of the right inverse property). Finally, every left deviation~\eqref{eq:padic-mobius-dev} lies in the unitary congruence group $U_2 := \{u\in 1+p^2\Oint : u\overline{u}=1\}$, which acts on $p\Oint/p^k\Oint$ with kernel $U_{k-1} := \{u\in 1+p^{k-1}\Oint : u\overline{u}=1\}$ (since $(u-1)p\equiv 0\pmod{p^k\Oint} \iff u\in 1+p^{k-1}\Oint$). As sets, the Cayley bijection $t\mapsto (1+it)(1-it)^{-1}$ maps $p^r\Zp$ onto $U_r$ for every $r\ge 1$, whence $|U_2/U_{k-1}| = p^{k-3}$ for all $k\ge 3$. Taking $\xi_1=p$ and $\xi_2=pi$, the factor~\eqref{eq:padic-mobius-dev} yields $u = (1+\varepsilon p^2 i)(1-\varepsilon p^2 i)^{-1} \equiv 1 + 2\varepsilon p^2 i \pmod{p^4\Oint}$. Since $p\ge 3$, binomial induction gives $u^{p^j} \equiv 1 + 2\varepsilon p^{2+j}i \pmod{p^{3+j}\Oint}$ for every $j\ge 0$, so $u\bmod U_{k-1}$ has exact order $p^{k-3}$ and generates $U_2/U_{k-1}$. Therefore $\Innl(L_k^{(\varepsilon)}) \cong U_2/U_{k-1} \cong \Z/p^{k-3}\Z$ is cyclic of exact order $p^{k-3}$ for all $k\ge 3$.
\end{proof}

\subsection{Higher \texorpdfstring{$p$-adic spheres ($n\ge 4$)}{p-adic spheres (n >= 4)}}

To achieve finite $A_\ell$-loops with \emph{non-abelian} left inner mapping groups of arbitrarily high nilpotency class, we pass to higher-dimensional $p$-adic spheres.

Let $p\ge 3$ be any odd prime and let $n\ge 4$. Set $R := M_n(\Zp)$ with the type~(I) conjugation involution $\sigma(X) := JXJ$ where $J := \mathrm{diag}(-1,1,\dots,1)$, and take the compact pro-$p$ congruence subgroup
\[
  G_1 := \SO_n(\Zp)_1 := \bigl\{ A\in I + pM_n(\Zp) : A^T A = I \bigr\}.
\]
(Note that $\det A = 1$ holds automatically for $A\in I+pM_n(\Zp)$ with $A^TA=I$, since $(\det A)^2 = 1$ and $\det A\equiv 1\pmod p$.) The Cayley transform $\kappa\colon p\,\so_n(\Zp)\xrightarrow{\sim} \SO_n(\Zp)_1$ is a global $p$-adic analytic diffeomorphism mapping $p^k\so_n(\Zp)$ onto $G_k := \SO_n(\Zp)_k$ for every $k\ge 1$. Under $\theta(X)=JXJ$, the Lie algebra $\mathcal{L} = \so_n(\Zp)$ splits as $\mathcal{L} = \h\oplus\m$, where
\[
  \h = \left\{ \begin{pmatrix} 0 & 0 \\ 0 & B \end{pmatrix} : B\in\so_{n-1}(\Zp) \right\}, \qquad
  \m = \left\{ \begin{pmatrix} 0 & v^T \\ -v & 0 \end{pmatrix} : v\in\Zp^{n-1} \right\}.
\]
Hence $H_1 = G_1^\tau \cong \SO_{n-1}(\Zp)_1$, $\rk_{\Zp}(\m) = n-1$, and $L = \kappa(p\m) = \{A\in\SO_n(\Zp)_1 : JAJ = A^{-1}\}$.

\begin{theorem}[Higher $p$-adic spheres and exact nilpotency class]\label{thm:higher-spheres}
Let $p\ge 3$ be any odd prime, let $n\ge 4$, and let $L_k := L/(L\cap\SO_n(\Zp)_k)$ ($k\ge 2$) be the finite Bruck/$A_\ell$-loop attached to $\SO_n(\Zp)_1/\SO_{n-1}(\Zp)_1$. Then:
\begin{enumerate}
  \item $|L_k| = p^{(n-1)(k-1)}$ for every $k\ge 2$;
  \item the faithful left inner mapping group $\Innl(L_k)$ has exact nilpotency class
  \[
    \cl\bigl(\Innl(L_k)\bigr) = \left\lfloor \frac{k-2}{2} \right\rfloor \qquad (\forall k\ge 2).
  \]
  In particular, $L_k$ is a proper Bruck/$A_\ell$-loop for all $k\ge 4$, has a non-abelian left inner mapping group $\Innl(L_k)$ for all $k\ge 6$, and its inverse limit $L=\varprojlim_k L_k$ has profinite left inner mapping group $\overline{\Innl(L)} \cong \varprojlim_k \Innl(L_k) \cong \SO_{n-1}(\Zp)_2$.
\end{enumerate}
\end{theorem}

\begin{proof}
Part~(1) follows from~\eqref{eq:finite-order} since $\rk_{\Zp}(\m)=n-1$. For~(2), write $M(v) := \bigl(\begin{smallmatrix} 0 & v^T \\ -v & 0\end{smallmatrix}\bigr)\in\m$ ($v\in\Zp^{n-1}$). A direct block computation yields
\[
  \bigl[M(u), M(v)\bigr] = \begin{pmatrix} 0 & 0 \\ 0 & vu^T - uv^T \end{pmatrix}.
\]
Define $E_{ij} := e_i e_j^T - e_j e_i^T$ for all distinct $i,j\in\{1,\dots,n-1\}$ (so $E_{ji}=-E_{ij}$). As $u,v$ range over the standard basis $e_1,\dots,e_{n-1}$ of $\Fp^{n-1}$, the skew-symmetric matrices $E_{ij}$ ($1\le i<j\le n-1$, arising from $[M(e_j),M(e_i)]$) form a basis of $\so_{n-1}(\Fp)$, whence $\bar{\h}_0 = [\bar{\m},\bar{\m}] = \so_{n-1}(\Fp)$. Since $n-1\ge 3$, for any distinct indices $i,j,r\in\{1,\dots,n-1\}$ one has $[E_{ir}, E_{rj}] = E_{ij}$ in $\so_{n-1}(\Fp)$ (valid over $\Fp$ for every odd prime $p\ge 3$, including $p=3$). Thus $[\bar{\h}_0,\bar{\h}_0] = \bar{\h}_0 \ne 0$, i.e., $\bar{\h}_0$ is perfect. Moreover, $\mathrm{ad}\bigl(\begin{smallmatrix}0&0\\0&B\end{smallmatrix}\bigr)(M(v)) = M(Bv)$, which is the standard faithful vector representation of $\so_{n-1}(\Fp)$ on $\Fp^{n-1}$. All hypotheses of Theorem~\ref{thm:exact-thresholds}(4) are therefore satisfied, yielding $\cl(\Innl(L_k)) = \lfloor(k-2)/2\rfloor$ for every $k\ge 2$.

Finally, for the inverse limit $L$, observe first that $H_2 = \SO_{n-1}(\Zp)_2$ acts faithfully by conjugation on $L = \kappa(p\m)$ (since $B\in\SO_{n-1}(\Zp)_2$ centralizes $\kappa(pM(v))$ if and only if $Bv=v$, whence $Bv=v$ for all $v\in\Zp^{n-1}$ forces $B=I_{n-1}$). Furthermore, since $p\ge 3$, for every $r\ge 2$ and every $B\in\so_{n-1}(\Zp)$ one has $(I+p^rB+O(p^{r+1}))^p \equiv I + p^{r+1}B \pmod{p^{r+2}}$ (as $v_p(\binom{p}{2}p^{2r}) = 2r+1 \ge r+3$), so $H_r^p$ surjects onto $H_{r+1}/H_{r+2}$. By induction and completeness of the closed subgroup $\overline{\langle H_2^p\rangle}$ (indeed $H_2$ is a uniform pro-$p$ group), it follows that $H_3 = \overline{\langle H_2^p\rangle} \le \Phi(H_2)$; since also $H_2^p\le H_3$ and $[H_2,H_2]\le H_4\le H_3$, the Frattini subgroup of $H_2$ is precisely $\Phi(H_2) = \overline{H_2^p[H_2,H_2]} = H_3$. Because $\overline{D_1}\le H_2$ is a closed subgroup mapping onto $H_2/H_3 \cong \so_{n-1}(\Fp) = \bar{\h}_0$, Burnside's basis theorem for pro-$p$ groups~\cite[Proposition~1.9]{dixon1999analytic} yields $\overline{D_1} = H_2$, whence $\overline{\Innl(L)} \cong \SO_{n-1}(\Zp)_2$.
\end{proof}

\subsection{\texorpdfstring{$p$-Adic polar decomposition ($n\ge 3$)}{p-Adic polar decomposition (n >= 3)}}

We conclude with a type~(II) example governed by the Cartan transpose-inverse involution, providing a $p$-adic analytic analogue of the classical real symmetric space $\mathrm{SL}_n(\R)/\mathrm{SO}_n(\R)$ of positive definite symmetric matrices~\cite{helgason1978differential}.

Let $p\ge 3$ be any odd prime and let $n\ge 3$. Set $R := M_n(\Zp)$ with the anti-automorphism $\sigma(X) := X^T$, so that on $G_1 := \GL_n(\Zp)_1 = I + pM_n(\Zp)$ the group involution is $\tau(A) := (A^T)^{-1}$ and $\theta(X) = -X^T$. The fixed subgroup is $H_1 = \OO_n(\Zp)_1 = \SO_n(\Zp)_1$, and the symmetric section is the closed subset of symmetric congruence matrices:
\[
  L := \bigl\{ P\in I + pM_n(\Zp) : P^T = P \bigr\} = \kappa\bigl(p\,\sym_n(\Zp)\bigr),
\]
where $\m = \sym_n(\Zp)$ has rank $\frac{n(n+1)}{2}$ and $\h = \so_n(\Zp)$. Theorem~\ref{thm:pro-p-cartan}(2) states that every matrix $A\in\GL_n(\Zp)_1$ admits a unique global \emph{$p$-adic polar decomposition} $A = P\cdot O$ with $P = (AA^T)^{1/2}\in L$ and $O = P^{-1}A\in\SO_n(\Zp)_1$, and the induced Bruck loop law on $L$ is $P\oplus Q = (PQ^2P)^{1/2}$.

Passing to determinant~$1$, notice that $\det\colon (L,\oplus)\to (1+p\Zp,\cdot)$ is a continuous loop epimorphism because $\det((PQ^2P)^{1/2}) = (\det(P)^2\det(Q)^2)^{1/2} = \det(P)\det(Q)$ (as $\det(P),\det(Q)\in 1+p\Zp$ and square roots are unique in the pro-$p$ group $1+p\Zp$). Consequently, the symmetric matrices of determinant~$1$,
\[
  L^{\SL} := L\cap\SL_n(\Zp)_1 = \ker(\det|_L),
\]
form a closed normal subloop of $(L,\oplus)$ (which coincides with the symmetric section of $\SL_n(\Zp)_1$ under $\tau$). Notice that although $\kappa(X)$ for traceless symmetric $X\in p\,\sym_n^0(\Zp)$ does not in general have determinant~$1$ when $n\ge 3$ (because $\det(I+X)\ne\det(I-X)$ at degree~$3$), the subgroup $G_1^{\SL} = \SL_n(\Zp)_1$ is closed and $\tau$-invariant. Moreover, for every $k\ge 2$, the finite quotient $L_k^{\SL} := L^{\SL}/(L^{\SL}\cap G_k)$ coincides with $\ker(\det_k\colon L_k\to 1+p(\Z/p^k\Z))$: indeed, if $P\in L$ satisfies $\det(P)\equiv 1\pmod{p^k}$, then $v := \det(P)^{-1/2}\in 1+p^k\Zp$ exists uniquely by Lemma~\ref{lem:unique-sqrt}, and setting $D_v := \mathrm{diag}(v,1,\dots,1)\in I+p^k M_n(\Zp)$ produces $P' := D_v P D_v \in L^{\SL}$ with $P'\equiv P\pmod{p^k M_n(\Zp)}$. Since $\det_k$ is surjective onto $1+p(\Z/p^k\Z)$ (already on $\mathrm{diag}(1+pt,1,\dots,1)$), it follows that $|L_k^{\SL}| = |L_k|/p^{k-1} = p^{(k-1)(\frac{n(n+1)}{2}-1)}$.

\begin{theorem}[$p$-Adic polar Bruck loops]\label{thm:polar-decomp}
Let $p\ge 3$ be any odd prime, let $n\ge 3$, and let $L_k := L/(L\cap\GL_n(\Zp)_k)$ and $L_k^{\SL} := L^{\SL}/(L^{\SL}\cap\SL_n(\Zp)_k)$ ($k\ge 2$). Then:
\begin{enumerate}
  \item $|L_k| = p^{(k-1)\frac{n(n+1)}{2}}$ and $|L_k^{\SL}| = p^{(k-1)\left(\frac{n(n+1)}{2}-1\right)}$ for all $k\ge 2$.
  \item For every $k\ge 2$, both $L_k$ and $L_k^{\SL}$ have exact left inner mapping nilpotency class
  \[
    \cl\bigl(\Innl(L_k)\bigr) = \cl\bigl(\Innl(L_k^{\SL})\bigr) = \left\lfloor \frac{k-2}{2} \right\rfloor.
  \]
\end{enumerate}
\end{theorem}

\begin{proof}
Part~(1) was established above via~\eqref{eq:finite-order} and the exact sequence $1\to L_k^{\SL}\to L_k\xrightarrow{\det_k} 1+p(\Z/p^k\Z)\to 1$. For~(2), in $\GL_n(\Zp)_1$ we have $\bar{\m} = \sym_n(\Fp)$ and $\bar{\h} = \so_n(\Fp)$. Taking the elementary symmetric matrices $S_{ij} := e_ie_j^T + e_je_i^T$ ($1\le i<j\le n$) and $D_{ij} := e_ie_i^T - e_je_j^T$ (all of which are traceless, hence lie in $\bar{\m}_0 := \sym_n^0(\Fp)$), we compute
\begin{align*}
  [D_{ij}, S_{ij}]
    &= (e_ie_i^T - e_je_j^T)(e_ie_j^T + e_je_i^T) - (e_ie_j^T + e_je_i^T)(e_ie_i^T - e_je_j^T) \\
    &= 2(e_ie_j^T - e_je_i^T) = 2E_{ij}.
\end{align*}
Since $2\in\Fp^\times$, the commutator subspace $[\bar{\m}_0,\bar{\m}_0]$ of traceless symmetric matrices already spans all of $\so_n(\Fp)$. Because $n\ge 3$ and $p\ge 3$, $\bar{\h}_0 = \so_n(\Fp)$ is perfect ($[E_{ir},E_{rj}]=E_{ij}$) and acts faithfully on $\bar{\m}_0 = \sym_n^0(\Fp)$: indeed, for any nonzero $B=(b_{ij})\in\so_n(\Fp)$ with $b_{ij}\ne 0$, the $(i,j)$-entry of $[B, D_{ij}]$ is $-2b_{ij}\ne 0$. Theorem~\ref{thm:exact-thresholds}(4) immediately yields $\cl(\Innl(L_k)) = \lfloor(k-2)/2\rfloor$.

For $L_k^{\SL}$, note first that $L^{\SL}\subseteq L$ inherits the upper bound $\cl(\Innl(L_k^{\SL}))\le \lfloor(k-2)/2\rfloor$ from $G_1$. Conversely, for every traceless symmetric matrix $X\in\sym_n^0(\Zp)$, expanding the characteristic polynomials $\det(I\pm pX) = 1 + p^2 e_2(X) \pm p^3 e_3(X) + O(p^4)$ (where $e_1(X)=\tr(X)=0$) shows that
\[
  \det\bigl(\kappa(pX)\bigr) = \frac{\det(I+pX)}{\det(I-pX)} \in 1 + p^3\Zp.
\]
Letting $u := \det(\kappa(pX))^{-1}\in 1+p^3\Zp$, Lemma~\ref{lem:unique-sqrt} provides a unique square root $v := u^{1/2}\in 1+p^3\Zp$. Setting $D_v := \mathrm{diag}(v,1,\dots,1)\in I + p^3 M_n(\Zp)$ and defining
\[
  P(X) := D_v\,\kappa(pX)\,D_v,
\]
one verifies at once that $P(X)^T = D_v\kappa(pX)^T D_v = P(X)$, $\det(P(X)) = v^2\det(\kappa(pX)) = 1$ (so $P(X)\in L^{\SL}$), and $P(X) \equiv \kappa(pX) \pmod{p^3 M_n(\Zp)}$ uniformly for every odd prime $p\ge 3$ (including when $p\mid n$). Since taking square roots in $1+pM_n(\Zp)$ preserves congruence modulo $p^3$, one has $P(X)\oplus P(Y) \equiv \kappa(pX)\oplus\kappa(pY) \pmod{p^3 M_n(\Zp)}$, and therefore
\[
  \delta\bigl(P(X),P(Y)\bigr) \equiv \delta\bigl(\kappa(pX),\kappa(pY)\bigr) \equiv I + 2p^2[X,Y] \pmod{p^3 M_n(\Zp)}.
\]
Thus the deviation subgroup $D_1^{\SL}\le H_2$ maps onto $\bar{\h}_0 = [\bar{\m}_0,\bar{\m}_0] = \so_n(\Fp)$ modulo $p^3$, and the induction of Theorem~\ref{thm:exact-thresholds}(4) applies verbatim to $L_k^{\SL}$ (using $P(Z)\equiv I+2pZ\pmod{p^2 M_n(\Zp)}$ for $Z\in\sym_n^0(\Zp)$), establishing $\cl(\Innl(L_k^{\SL})) = \lfloor(k-2)/2\rfloor$ for every $k\ge 2$ and every odd prime $p\ge 3$.
\end{proof}

\begin{remark}\label{rmk:padic-summary}
Table~\ref{tab:padic-examples} collects the exact invariants of our three families at the smallest odd prime $p=3$ across congruence levels $k=2,\dots,7$ (all entries have been checked independently by modular matrix computation over $\Z/3^k\Z$ up to $k=8$, where $\cl(\Innl(L_8))=3$ in the last three rows). Together with the left-folder correspondence of Baer and Sabinin~\cite{baer1939nets,sabinin1972equivalence}, these $p$-adic symmetric towers provide a systematic source of finite $A_\ell$-loops whose left inner mapping groups attain arbitrary nilpotency class, complementing constructions derived from split semidirect products.

\begin{table}[ht]
\centering
\small
\caption{Invariants of finite $p$-adic Bruck/$A_\ell$-loops at $p=3$.}\label{tab:padic-examples}
\begin{tabular}{@{}lcccc@{}}
\toprule
Family & $|L_k|$ & $k=2,3$ & $k=4,5$ & $k=6,7$ \\
\midrule
M\"obius ($p\Z_3[i]/3^k$) & $3^{2(k-1)}$ & Abelian & $\Innl\cong\Z/3^{k-3}\Z$ & $\Innl\cong\Z/3^{k-3}\Z$ \\
$\SO_4(\Z_3)_1/\SO_3(\Z_3)_1$ & $3^{3(k-1)}$ & $\cl=0$ & $\cl=1$ (abelian) & $\cl=2$ (non-abelian) \\
$\GL_3(\Z_3)_1/\OO_3(\Z_3)_1$ & $3^{6(k-1)}$ & $\cl=0$ & $\cl=1$ (abelian) & $\cl=2$ (non-abelian) \\
$\SL_3(\Z_3)_1/\SO_3(\Z_3)_1$ & $3^{5(k-1)}$ & $\cl=0$ & $\cl=1$ (abelian) & $\cl=2$ (non-abelian) \\
\bottomrule
\end{tabular}
\end{table}
\end{remark}

\bibliographystyle{amsplain}
\bibliography{references_beta}

\end{document}